\documentclass[12pt,a4paper]{amsart}
\usepackage[utf8]{inputenc}	
\usepackage[T1]{fontenc}
\usepackage[in,headings]{fullpage}
\usepackage{amsbsy, amscd, amsfonts, amsmath, amsrefs, amssymb, amsthm}
\usepackage{graphicx}
\usepackage{float}
\usepackage{color}   %May be necessary if you want to color links
\usepackage[colorlinks=true,linkcolor=blue,citecolor=blue,urlcolor=blue]{hyperref}
\usepackage{comment}
\excludecomment{A}

\newtheorem{theorem}{Theorem}%[section]

\newtheorem{lemma}[theorem]{Lemma}

\theoremstyle{definition}
\newtheorem{definition}[theorem]{Definition}

\theoremstyle{remark}

\newcommand{\C}{\mathbf{C}}
\newcommand{\Z}{\mathbf{Z}}

\newcommand{\R}{\mathbf{R}}
\newcommand{\N}{\mathbf{N}}

\renewcommand{\Re}{\mathop{\mathrm{Re}}\nolimits}
\renewcommand{\Im}{\mathop{\mathrm{Im}}\nolimits}
\newcommand{\Rzeta}{\mathop{\mathcal R }\nolimits}

\newfont{\cmbsy}{cmbsy10}
\newfont{\cmmib}{cmmib10}
\newcommand{\Orden}{\mathop{\hbox{\cmbsy O}}\nolimits}

\DeclareMathOperator{\sgn}{sgn}
\DeclareMathOperator{\Res}{Res}

\begin{document}

\title[Approximate formula for $Z(t)$]
{Approximate formula for $Z(t)$}
\author[Arias de Reyna]{J. Arias de Reyna}
\address{%
Universidad de Sevilla \\ 
Facultad de Matem\'aticas \\ 
c/Tarfia, sn \\ 
41012-Sevilla \\ 
Spain.} 

% AMS subject classifications (used in AMS journals)

% AMS keywords (used in AMS journals)

% acknowledge support, etc
%\thanks{This research was  supported by MINECO grant MTM2015--63699-P}
% \thanks{We would like to thank our colleagues for their helpful
%  criticism.}

\email{arias@us.es, ariasdereyna1947@gmail.com}

%\date{\today, \texttt{60-Approximate-v3.tex}}

\begin{abstract}
The series for the zeta function does not converge on the critical line but the function \[G(t)=\sum_{n=1}^\infty \frac{1}{n^{\frac12+it}}\frac{t}{2\pi n^2+t},\] satisfies $Z(t)=2\Re\{e^{i\vartheta(t)}G(t)\}+\Orden(t^{-\frac56+\varepsilon})$. So one expects that the zeros of zeta on the critical line are very near the zeros of $\Re\{e^{i\vartheta(t)}G(t)\}$. There is a related function $U(t)$ that satisfies the equality $Z(t)=2\Re\{e^{i\vartheta(t)}U(t)\}$. 
\end{abstract}

\maketitle

\setcounter{section}{0}
\section*{Introduction.}

There are many functions $f(t)$ with $Z(t)=2\Re\{e^{i\vartheta(t)}f(t)\}$. Notable examples of this type of functions are $\zeta(\frac12+it)$, $-\frac{\zeta'(\frac12+it)}{2\vartheta'(t)}$, and $\Rzeta(\frac12+it)$. The position of the zeros of these functions is closely related to the zeros of zeta on the critical line. In this paper we define a function $G(t)$ that only satisfies the approximate equation
\[Z(t)=2\Re\{e^{i\vartheta(t)}G(t)\}+\Orden(t^{-\frac56+\varepsilon}).\]
Therefore, we expect that the real zeros of $2\Re\{e^{i\vartheta(t)}G(t)\}$ are close to the real zeros of $Z(t)$. In exchange for the fact that it is only an approximation, we find that $G(t)$ has a very simple definition, valid for all real $t$, \[G(t)=\sum_{n=1}^\infty \frac{1}{n^{\frac12+it}}\frac{t}{2\pi n^2+t},\] and that $2\Re\{e^{i\vartheta(t)}G(t)\}$ seems to have of order of $\frac{T}{2\pi}\log T$ real zeros in the interval $[0,T]$. Although the latter we have not been able to prove, I only base this statement on the calculation of some of the first zeros of $G(t)$. 

The function $G(t)$ extends to a meromorphic function on $\C$ with poles at $(2k-\frac12)i$ and $-2\pi n^2$. There are many interesting things about $G(t)$ for example
\begin{equation}
G(t)=\frac{t}{2\pi}\Bigl(\zeta(2+\tfrac12+it)-\sum_{n=1}^\infty \frac{1}{n^{2+\frac12+it}}
\frac{t}{2\pi n^2+t}\Bigr).
\end{equation}
So that the values of $Z(t)$ appear to depend only on things happening at $\sigma\ge2$.
Also for small $t$ we observe that the continuous  $\arg G(t)$ is increasing very slowly. Giving the zeros we want for $Z(t)$.
\begin{figure}[H]
\begin{center}
\includegraphics[width=10cm]{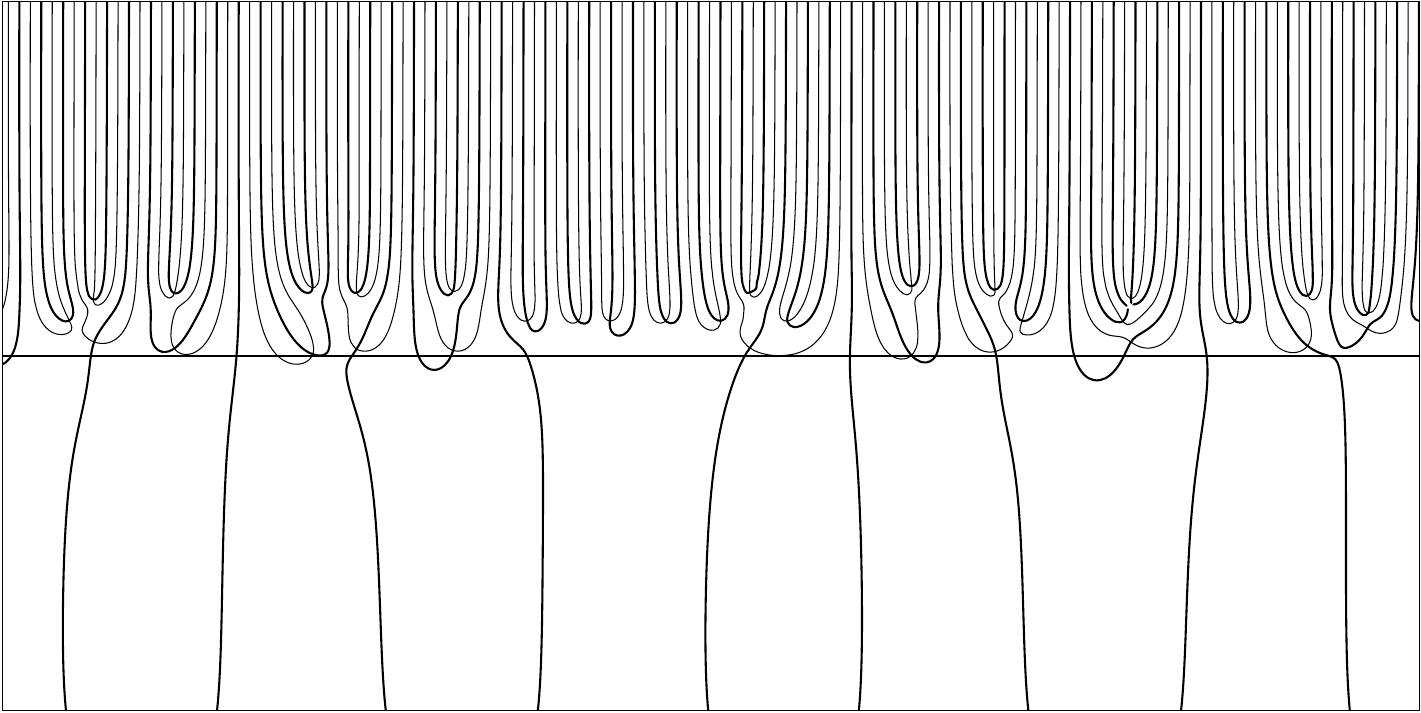}
\caption{Plot of  $G(t)$ in $(1000,1040)\times(-10,10)$.}
\label{secondplot4}
\end{center}
\end{figure}

The function $G(t)$ is only a first example of many other approximate solutions. For example, see \cite{A162}.

\subsection*{Notations and known results.} 
Throughout the paper, $C$'s and $c$'s denote positive absolute constants that are not always the same from one occurrence to another. $\Orden$'s are absolute unless otherwise stated. Thus $f(x)=\Orden(\phi(x))$ means $|f(x)|/\phi(x)< K$ for $x>x_0$: an $\Orden$ involves two constants, $K$ and $x_0$. We also write $f(x)\ll \phi(x)$ with the same meaning.

We recall some known results. First,  by Stirling's formula (see Titchmarsh \cite{T}*{(4.12.1), p.~78}. In any fixed strip
$\alpha\le\sigma\le\beta$, as $t\to\infty$
\begin{equation}\label{Stirling}
\Gamma(\sigma+it)=(2\pi)^{1/2} |t|^{\sigma+it-1/2}e^{-\pi
|t|/2-it+\sgn(t)\frac{\pi
i}{2}(\sigma-1/2)}\bigl(1+\Orden(|t|^{-1})\bigr).
\end{equation}

The function $Z(t)$ is defined (see \cite{T}*{Section 4.17}, \cite{Edwards}*{Section 6.5}) as
\begin{equation}\label{defZ}
Z(t)=e^{i\vartheta(t)}\zeta(\tfrac12+it),\quad \text{ where }\quad
e^{i\vartheta(t)}=\pi^{-it/2}\sqrt{\frac{\Gamma\bigl(\frac{1}{4}+i\frac{t}{2}\bigr)}
{\Gamma\bigl(\frac{1}{4}-i\frac{t}{2}\bigr)}}=\frac{\pi^{-\frac12it}\Gamma(\frac14+\frac12it)}{|\Gamma(\frac14+\frac12it)|}.
\end{equation}
It is an analytical function on $\Omega$  the plane minus two cuts along 
the imaginary axis, one from $i/2$ to $i\infty$ and one from $-i/2$
to $-i\infty$.  The
function $\vartheta(t)$ is also holomorphic on $\Omega$. It is
defined so that $\vartheta(0)=0$.  The functional equation of zeta implies 
that $Z(t)$ is an even function and $\vartheta(t)$ an odd function.

By the Stirling formula, for  $t=x+iy$ and $\alpha\le y\le \beta$ and
$x\to\infty$ we have
\begin{equation}\label{BoundVartheta}
e^{i\vartheta(x+iy)}=\Bigl|\frac{x}{2\pi}\Bigr|^{-y/2+ix/2}e^{-i\frac{x}{2}
-\frac{\pi i}{8}\sgn(x)}\Bigl(1+\Orden (|x|^{-1})\Bigr).
\end{equation}
We also need a simple bound for $\zeta(\sigma+it)$ Edwards \cite{Edwards}*{p.~185}
\begin{equation}\label{simpleboundzeta}
|\zeta(\sigma+it)|\le Ct^{\frac{1-\sigma}{2}}\log t,\qquad 0\le\sigma\le 1,\quad t>2.
\end{equation}
Sometimes, it suffices to apply the more simple inequality
\begin{equation}\label{boundzeta}
|\zeta(\sigma+it)|\le Ct^{\frac12},\qquad \sigma\ge \tfrac12,\quad t>1,
\end{equation}
that we may deduce from \cite{T}*{Theorem 4.11}.

For the incomplete Gamma function we have the following inequality (see Gabcke 
\cite{G}*{Ch.~5 \S~4 Theorem 3})
\begin{equation}\label{GammaInc}
\Gamma(a,x):=\int_x^\infty v^{a-1}e^{-v}\,dv\le ae^{-x}x^{a-1},\qquad x>a\ge1.
\end{equation}

We shall make use of the following asymptotic expansions  (see Gabcke \cite{G}*{Ch.~4 \S~2 Theorem 3})
\begin{equation}\label{AsympTheta}
\vartheta(t)\sim\frac{t}{2}\log\frac{t}{2\pi}-\frac{t}{2}-\frac{\pi}{8}
+\sum_{n=1}^\infty\frac{(2^{2n-1}-1)|B_{2n}|}{2^{2n}(2n-1)
2n}\frac{1}{t^{2n-1}}.
\end{equation}
\begin{equation}\label{AsympThetaPrime}
\vartheta'(t)\sim\frac{1}{2}\log\frac{t}{2\pi}
-\sum_{n=1}^\infty\frac{(2^{2n-1}-1)|B_{2n}|}{2^{2n}
2n}\frac{1}{t^{2n}}.
\end{equation}
\begin{equation}\label{AsympThetaSegunda}
\vartheta''(t)\sim\frac{1}{2t}+
\sum_{n=1}^\infty\frac{(2^{2n-1}-1)|B_{2n}|}{2^{2n}
}\frac{1}{t^{2n+1}}.
\end{equation}
Since they follow from the Stirling expansion, all of them are valid on $|\arg t|<\theta<\frac{\pi}{2}$.

\section{The function \texorpdfstring{$U(t)$}{U(t)}.}

\begin{definition}
For $t\in\R$ we define
\begin{equation}\label{defU}
U(t)=\frac{1}{2\pi i}\int\limits_{-i\sigma-\infty}^{-i\sigma+\infty}
e^{i\vartheta(t+x)-i\vartheta(t)}\zeta\bigl(\tfrac12+i(t+x)\bigr)\frac{\pi}{2\sinh\frac{\pi x}{2}}\,dx,
\end{equation}
where $0<\sigma<\frac12$.
\end{definition}

It is easy to see that the integrand is a holomorphic function of  $x$ for 
$-\frac12<\Im(x)<0$, and by Cauchy's Theorem the integral is independent of $\sigma$. 

\begin{theorem}\label{T:ZU}
For $t\in\R$ we have
\begin{equation}
Z(t)=2\Re\bigl\{e^{i\vartheta(t)}U(t)\bigr\}.
\end{equation}
\end{theorem}

\begin{proof}
Following an idea of Berry and Keating \cite{BK}, by Cauchy's Theorem, we have
\begin{displaymath}
Z(t)=\frac{1}{2\pi i} \int_{L_-+L_+}Z(t+x)\frac{\pi}{2\sinh\frac{\pi x}{2}}\,dx,
\end{displaymath}
where for $0<\sigma<\frac12$, the integration path  $L^-$ and $L^+$ are  parametrized by $x=-i\sigma+iy$ and $x= i\sigma-iy$ respectively with $-\infty<y<+\infty$. 

That is,
\begin{displaymath}
Z(t)=\frac{1}{2\pi i} \int\limits_{-i\sigma-\infty}^{-i\sigma+\infty}Z(t+x)\frac{\pi}{2\sinh\frac{\pi x}{2}}\,dx-\frac{1}{2\pi i}\int\limits_{i\sigma-\infty}^{i\sigma+\infty}Z(t+x)\frac{\pi}{2\sinh\frac{\pi x}{2}}\,dx.
\end{displaymath}
For real $t$,  the second integral is the conjugate of the first integral, therefore
\begin{displaymath}
Z(t)=2\Re\Bigl\{\frac{1}{2\pi i} \int\limits_{-i\sigma-\infty}^{-i\sigma+\infty}Z(t+x)\frac{\pi}{2\sinh\frac{\pi x}{2}}\,dx\Bigr\}.
\end{displaymath}
This may be written as 
\begin{displaymath}
Z(t)=2\Re\Bigl\{e^{i\vartheta(t)}\frac{1}{2\pi i} \int\limits_{-i\sigma-\infty}^{-i\sigma
+\infty}e^{i\vartheta(t+x)-i\vartheta(t)}\zeta\bigl(\tfrac12+i(t+x)\bigr)
\frac{\pi}{2\sinh\frac{\pi x}{2}}\,dx\Bigr\}.\qedhere
\end{displaymath}

\end{proof}

\section{First approximation to \texorpdfstring{$U(t)$}{U(t)}.}

\begin{theorem}\label{firstT}
Given a natural number $n$ and $0<\sigma_0<2$ there is a constant $C>0$ such 
that for $t\to+\infty$
\begin{equation}
U(t)=\frac{1}{2\pi i}\int\limits_{-b-i\sigma_0}^{b-i\sigma_0}e^{i\vartheta(t+x)-i\vartheta(t)}\zeta\bigl(\tfrac12+i(t+x)\bigr)\frac{\pi}{2\sinh\frac{\pi x}{2}}\,dx+\Orden(t^{-n}),
\end{equation}
with $b=C\log t$. 
\end{theorem}

\begin{proof}
In equation \eqref{defU} we may take $0<\sigma<\tfrac12$ and $\sigma<\sigma_0$.  Applying 
Cauchy's Theorem we change the path of integration to a broken line $L$ through the points
\begin{displaymath}
-\infty-i\sigma, \quad -b-i\sigma,\quad -b-i\sigma_0,\quad +b-i\sigma_0,\quad 
+\infty-i\sigma_0.
\end{displaymath}
We call $L_1$, $L_2$, $L_3$, and $L_4$ the four segments that make up $L$. 
Since we want to approximate $U(t)$ by the integral on $L_3$, we must bound the other 
three integrals.
\medskip

(a) Integral along $L_4$. 
The path is parametrized by $x=y-i\sigma_0 $ with $b<y<+\infty$. 
By \eqref{BoundVartheta}
\begin{displaymath}
|e^{i\vartheta(t+x)}|=|e^{i\vartheta(t+y-i\sigma_0)}|\ll \Bigl|\frac{t+y}{2\pi}\Bigr|
^{\frac{\sigma_0}{2}}\ll|t+y|^{\frac{\sigma_0}{2}},
\end{displaymath}
and by \eqref{boundzeta} since $y>b=C\log t>1$ 
\begin{displaymath}
|\zeta\bigl(\tfrac12+i(t+x)\bigr)|=|\zeta\bigl(\tfrac12+\sigma_0+i(t+y)\bigr)|
\le C(t+y)^{\frac{1}{2}}.
\end{displaymath}
Since $0<\sigma_0<2$ the last factor is bounded by
\begin{displaymath}
\Bigl|\frac{\pi}{2\sinh\frac{\pi (y-i\sigma_0)}{2}}\Bigr|\ll e^{-\frac{\pi y}{2}}.
\end{displaymath}
It follows that 
\begin{multline*}
I(L_4):=\Bigl|\frac{1}{2\pi i}\int\limits_{b-i\sigma_0}^{\infty-i\sigma_0}e^{i\vartheta(t+x)
-i\vartheta(t)}\zeta\bigl(\tfrac+i(t+x)\bigr)\frac{\pi}{2\sinh\frac{\pi x}{2}}
\,dx\Bigr|\ll\\
\ll \int_b^{+\infty}(t+y)^{\frac{1+\sigma_0}{2}}e^{-\frac{\pi y}{2}}\,dy
\ll \int_b^{+\infty}(t+y)^{2} e^{-\frac{\pi y}{2}}\,dy.
\end{multline*}
Since $1<b<t$ we have $t+y\le ty$ and 
\begin{displaymath}
I(L_4)\ll t^{2} \int_b^{+\infty}y^{2} e^{-\frac{\pi y}{2}}\,dy
\ll  t^{2} \int_{\pi b/2}^{+\infty}y^{2} e^{-y}\,dy
= t^{\frac12}\Gamma(3,\pi b/2).
\end{displaymath}
Now we use \eqref{GammaInc} and $b=C\log t$
\begin{displaymath}
I(L_4)\ll t^{2}e^{-\frac{\pi b}{2}}\Bigl(\frac{\pi b}{2}\Bigr)^{2}
\ll t^{2}t^{-\frac{\pi C}{2}}\Bigl(\frac{\pi C}{2}\log t\Bigr)^{2}\ll t^{-n},
\end{displaymath}
if we take $C=cn$ with $c$ large enough.

(b) Integral along $L_2$.  The path  is parametrized by $x=-b-iy$ with $\sigma<y<\sigma_0$.

By \eqref{BoundVartheta} we have the following
\begin{displaymath}
|e^{i\vartheta(t+x)}|=|e^{i\vartheta(t-b-iy)}|\ll\Bigl|\frac{t-b}{2\pi}\Bigr|^{y/2}.
\end{displaymath}
We choose $b=C\log t$ and $t$ large enough so that $t-b>2\pi$. Also, observe that
$y<\sigma_0<2$ and  hence
\begin{equation}
|e^{i\vartheta(t+x)}|\ll (t-b)<t.
\end{equation}
We have
\begin{equation}
|\zeta\bigl(\tfrac12+i(t+x)\bigr)|=|\zeta\bigl(\tfrac12+y+i(t-b)\bigr)|\ll t^{\frac12}.
\end{equation}
On the other hand,
\begin{equation}
\Bigl|\frac{\pi}{2\sinh\frac{\pi x}{2}}\Bigr|=\Bigl|\frac{\pi}
{2\sinh\frac{\pi (b+iy)}{2}}\Bigr|\ll e^{-\frac{\pi b}{2}}.
\end{equation}
Hence,
\begin{equation}
I(L_2)\ll \int_{\sigma}^{\sigma_0} t^{\frac32}t^{-\frac{\pi C}{2}}\,dy\ll
t^{\frac32}t^{-\frac{\pi C}{2}}=\Orden(t^{-n}).
\end{equation}
\medskip

(c) Integral along $L_1$. The path is parametrized by $x=y-i\sigma$ with $-\infty <y<-b$. In this case, we have 
\begin{displaymath}
|e^{i\vartheta(t+x)}|=|e^{i\vartheta(t+y-i\sigma)}|,
\end{displaymath}
but now $t+y$ ranges through $(-\infty, t-b)$ so it takes small values. Since $0<\sigma<\frac12$ the function $\vartheta(t+y-i\sigma)$ is continuous and for $|t+y|$ big it is 
given by \eqref{BoundVartheta}. It follows that 
\begin{equation}
|e^{i\vartheta(t+x)}|=|e^{i\vartheta(t+y-i\sigma)}|\ll (1+|t+y|^{\frac{\sigma}{2}}).
\end{equation}
In the same way 
\begin{equation}
|\zeta\bigl(\tfrac12+i(t+x)\bigr)|=|\zeta\bigl(\tfrac12+\sigma+i(t+y)\bigr)|\ll (1+|t+y|^{\frac12}).
\end{equation}
And
\begin{equation}
\Bigl|\frac{\pi}{2\sinh\frac{\pi x}{2}}\Bigr|=\Bigl|\frac{\pi}
{2\sinh\frac{\pi (y-i\sigma)}{2}}\Bigr|\ll e^{\frac{\pi y}{2}}.
\end{equation}
Hence,
\begin{multline*}
I(L_1)\ll \int_{-\infty}^{-b}(1+|t+y|^{\frac12})
(1+|t+y|^{\frac{\sigma}{2}})e^{\frac{\pi y}{2}}\,dy=\\
=\int_b^{\infty}
(1+|t-y|^{\frac12}+|t-y|^{\frac{\sigma}{2}}+|t-y|^{\frac{1+\sigma}{2}})
e^{-\frac{\pi y}{2}}\,dy.
\end{multline*}
For each one $|t-y|^\alpha\le 1+|t|+|y|$ since $|t-y|\le1$ or else 
$|t-y|^\alpha\le |t-y|$ for $\alpha\le 1$. Therefore,
\begin{align*}
I(L_1)&\ll\int_b^\infty(1+t+y)e^{-\frac{\pi y}{2}}\,dy=
\Bigl(\frac{4}{\pi^2}+\frac{2(1+b+t)}{\pi}\Bigr)e^{-\frac{\pi b}{2}}\\
&\ll(t+C\log t+1)t^{-\frac{\pi C}{2}}=\Orden(t^{-n}),
\end{align*}
by the choice of $C$. 
\end{proof}

For the next Theorem we shall need a lemma.
\begin{lemma}\label{L:uno}
For $t>0$ real, $|z|\le t/4$ and when  $t\to\infty$ we have
\begin{equation}
\vartheta(t+z)=\vartheta(t)+z\vartheta'(t)+\frac{z^2}{2}\vartheta''(t)+
\Orden\Bigl(\frac{\log |t|}{|t|^2}|z|^3\Bigr).
\end{equation}
\end{lemma}

\begin{proof}
By Cauchy's Theorem
\begin{displaymath}
\vartheta(t+z)-\vartheta(t)-z\vartheta'(t)-\frac{1}{2}\vartheta''(t)
z^2=\frac{z^3}{2\pi
i}\int_{C_R}\frac{\vartheta(t+\zeta)}{\zeta-z}\frac{d\zeta}{\zeta^3},
\end{displaymath}
when $C_R$ is the boundary of a disc that contains the two poles
$\zeta=0$ y $\zeta=z$ of the integrand. If we assume that
$|z|<t/4$ we can take  $C_R$ as the disc with center in $0$ and
radius  $|t|/2$.

At the points of the circle $|t|/2\le |t+\zeta|\le 3|t|/2$,  thus
by the asymptotic expansion  \eqref{AsympTheta} we get
$|\vartheta(t+\zeta)|\le C |t|\log|t|$. Hence we will have the following
\begin{displaymath}
\left|\vartheta(t+z)-\vartheta(t)-z\vartheta'(t)-\frac{1}{2}\vartheta''(t)
z^2\right|\le C \frac{|z|^3}{2\pi}\frac{|t|\log|t|}{|t|^4}|t|\le C
\log |t|\frac{|z|^3}{|t|^2}.\qedhere
\end{displaymath}
\end{proof}

\begin{theorem}\label{T:secodapr}
Let $0<\varepsilon<1$ and  $0<\sigma_0<2$,  then 
with the notation of Theorem \ref{firstT}, and taking $C>0$ big enough, we have
for $t\to+\infty$
\begin{equation}\label{secodapr}
U(t)=\frac{1}{2\pi i}\int\limits_{-b-i\sigma_0}^{b-i\sigma_0}e^{ix\vartheta'(t)+i\frac{x^2}{2}\vartheta''(t)}\zeta\bigl(\tfrac12+i(t+x)\bigr)\frac{\pi}{2\sinh\frac{\pi x}{2}}\,dx+\Orden(t^{-\frac74+\varepsilon}),
\end{equation}
with $b=C\log t$. 
\end{theorem}

\begin{proof}
First, we consider the case $\sigma_0=\frac12+\varepsilon$. 
Take $C$ large enough so that Theorem \ref{firstT} applies with $n=2$ and this $\sigma_0=\frac12+\varepsilon$. 
Hence, $L_3$ being the segment with extremes $-b-i\sigma_0$ and $b-i\sigma_0$, 
\begin{equation}\label{E:21}
U(t)=\frac{1}{2\pi i}\int_{L_3}e^{i\vartheta(t+x)-i\vartheta(t)}
\zeta\bigl(\tfrac12+i(t+x)\bigr)\frac{\pi}{2\sinh\frac{\pi x}{2}}\,dx+\Orden(t^{-2}).
\end{equation}
We call $I(t)$ the above integral and $I_0(t)$ that in \eqref{secodapr}. 
By Lemma \ref{L:uno} we have
\begin{gather}
\vartheta(t+x)-\vartheta(t)=x\vartheta'(t)+\frac{x^2}{2}\vartheta''(t)+w(t,x),\\
w(t,x)=\Orden\Bigl(\frac{\log t}{t^2}|x|^3\Bigr),\quad |x|<t/4.\label{wbound}
\end{gather}

For all $x$ in $L_3$, the  bound \eqref{wbound} applies. 
Therefore,
\begin{equation}
I(t)-I_0(t)=\int_{L_3}e^{ix\vartheta'(t)+i\frac{x^2}{2}\vartheta''(t)}(e^{iw(t,x)}-1)
\zeta\bigl(\tfrac12+i(t+x)\bigr)\frac{\pi}{2\sinh\frac{\pi x}{2}}\,dx.
\end{equation}
We need to bound this difference. To this end, observe that the path $L_3$ is parametrized by $x=y-i\sigma_0$ with $-b<y<b$.  So that by \eqref{AsympThetaPrime} and \eqref{AsympThetaSegunda}
\begin{equation}
|e^{ix\vartheta'(t)}|=e^{\sigma_0\vartheta'(t)}\ll \Bigl(\frac{t}{2\pi}\Bigr)^{\sigma_0/2}.
\end{equation}
\begin{equation}
|e^{i\frac{x^2}{2}\vartheta''(t)}|=e^{\sigma_0\vartheta''(t) y}\le 
e^{c\sigma_0\frac{|y|}{t}}\ll\Bigl(1+\frac{\log t}{t}\Bigr)\ll1.
\end{equation}
\begin{equation}
|e^{iw(t,x)}-1|\ll\frac{\log t}{t^2}y^3\ll\frac{\log^4t}{t^2}.
\end{equation}
We assume that $\sigma_0=\frac12+\varepsilon$ so that
\begin{equation}
|\zeta\bigl(\tfrac12+i(t+x)\bigr)|=|\zeta\bigl(\tfrac12+\sigma_0+i(t+y)\bigr)|\le 
\zeta\bigl(1+\varepsilon)\ll1.
\end{equation}
and
\begin{equation}
\Bigl|\frac{\pi}{2\sinh\frac{\pi x}{2}}\Bigr|=\Bigl|\frac{\pi}
{2\sinh\frac{\pi (y-i\sigma_0)}{2}}\Bigr|\ll e^{-\frac{\pi |y|}{2}}.
\end{equation}
Therefore,
\[
|I(t)-I_0(t)|\ll \int_{-b}^b  t^{\frac{\sigma_0}{2}}\frac{\log^4t}{t^2}e^{-\frac{\pi|y|}{2}}\,dy \ll t^{-(2-\frac{\sigma_0}{2})}\log^4t=t^{-\frac74+\frac{\varepsilon}{2}}
\log^4t\ll t^{-\frac74+\varepsilon}.\label{E:29}
\]

By \eqref{E:21} we have $U(t)=I(t)+\Orden(t^{-2})$   with \eqref{E:29} we have proved
$U(t)=I_0(t)+\Orden(t^{-\frac74+\varepsilon})$.  Therefore, we have \eqref{secodapr} for 
$\sigma_0=\frac12+\varepsilon$. 

Put $I_0(t)=I_0(t,\frac12+\varepsilon)$ to get our general result we only have
to show that for any $0<\sigma_1<\sigma_2<2$ we have $I_0(t,\sigma_1)=I_0(t,\sigma_2)+\Orden(t^{-2})$.  By Cauchy's Theorem we have
\begin{displaymath}
I_0(t,\sigma_2)-I_0(t,\sigma_1)=-J_0(t,b)+J_0(t,-b),
\end{displaymath}
where 
\begin{equation}
J_0(t,a):=\frac{1}{2\pi i}\int\limits_{a-i\sigma_2}^{a-i\sigma_1}
e^{ix\vartheta'(t)+i\frac{x^2}{2}\vartheta''(t)}
\zeta\bigl(\tfrac12+i(t+x)\bigr)\frac{\pi}{2\sinh\frac{\pi x}{2}}\,dx.
\end{equation}
However,  both integrals $J_0(t,b)$ and $J_0(t,-b)$ are $\Orden(t^{-2})$. The two proofs are similar. For example, for $J_0(t,b)$ the integration path is given by $x=b+iy$ with 
$-\sigma_2<y<-\sigma_1$. Hence $0<-y<2$ and we have
\begin{equation}
|e^{ix\vartheta'(t)}|=e^{-y\vartheta'(t)}\ll\Bigl(\frac{t}{2\pi}\Bigr)^{-y/2}\ll t.
\end{equation}
\begin{equation}
|e^{i\frac{x^2}{2}\vartheta''(t)}|=e^{-by\vartheta''(t)}\ll e^{2C\frac{\log t}{t}}\ll1.
\end{equation}
\begin{equation}
|\zeta\bigl(\tfrac12+i(t+x)\bigr)|=|\zeta\bigl(\tfrac12-y+i(t+b)\bigr)|\ll(t+b)^{\frac12}
\ll t^{\frac12}.
\end{equation}
and
\begin{equation}
\Bigl|\frac{\pi}{2\sinh\frac{\pi x}{2}}\Bigr|=\Bigl|\frac{\pi}
{2\sinh\frac{\pi (b+iy)}{2}}\Bigr|\ll e^{-\frac{\pi b}{2}},
\end{equation}
so that
\begin{equation}
|J_0(t,b)|\ll t^{\frac32}e^{-\frac{\pi b}{2}}\ll t^{\frac32} t^{-\frac{\pi C}{2}}
=\Orden(t^{-2}).\qedhere
\end{equation}
\end{proof}

\begin{theorem}\label{T:3}
Let $0<\varepsilon<1$ and  $0<\sigma<2$,  then 
with the notation of Theorem \ref{firstT}, and taking $C>0$ big enough, we have
for $t\to+\infty$
\begin{equation}\label{3apr}
U(t)=\frac{1}{2\pi i}\int\limits_{-b-i\sigma}^{b-i\sigma}
\Bigl(\frac{t}{2\pi}\Bigr)^{\frac{ix}{2}}
\zeta\bigl(\tfrac12+i(t+x)\bigr)\frac{\pi}{2\sinh\frac{\pi x}{2}}\,dx+\Orden(t^{-\frac56+\varepsilon}),
\end{equation}
with $b=C\log t$. 

Assuming Lindelöf hypothesis we may substitute the error term in \eqref{3apr} by $\Orden(t^{-1+\varepsilon})$.
\end{theorem}

\begin{proof}
Take $\sigma_0=\varepsilon$ in Theorem \ref{T:secodapr}. Therefore,
\begin{displaymath}
U(t)=I_0(t,\sigma_0)+\Orden(t^{-\frac74+\varepsilon}).
\end{displaymath}
where $I_0(t,\sigma_0)$ is the integral in \eqref{secodapr}.  Now, define
\begin{equation}
K(t,\sigma_0):=\frac{1}{2\pi i}\int\limits_{-b-i\sigma_0}^{b-i\sigma_0}
\Bigl(\frac{t}{2\pi}\Bigr)^{\frac{ix}{2}}
\zeta\bigl(\tfrac12+i(t+x)\bigr)\frac{\pi}{2\sinh\frac{\pi x}{2}}\,dx.
\end{equation}
Then 
\begin{multline*}
I_0(t,\sigma_0)-K(t,\sigma_0)=\frac{1}{2\pi i}
\int\limits_{-b-i\sigma_0}^{b-i\sigma_0}
\Bigl(\frac{t}{2\pi}\Bigr)^{\frac{ix}{2}}
\Bigl\{e^{ix\vartheta'(t)+i\frac{x^2}{2}\vartheta''(t)-i\frac{x}{2}\log\frac{t}{2\pi}}
-1\Bigr\}\cdot\\
\cdot
\zeta\bigl(\tfrac12+i(t+x)\bigr)\frac{\pi}{2\sinh\frac{\pi x}{2}}\,dx.
\end{multline*}
The path of integration is given by $x=y-i\sigma_0=y-i\varepsilon$ and therefore
\begin{equation}
\Bigl|\Bigl(\frac{t}{2\pi}\Bigr)^{\frac{ix}{2}}\Bigr|=
\Bigl(\frac{t}{2\pi}\Bigr)^{\frac{\sigma_0}{2}}\ll t^{\varepsilon/2}.
\end{equation}
By \eqref{AsympThetaPrime} and \eqref{AsympThetaSegunda}
\begin{displaymath}
ix\vartheta'(t)+i\frac{x^2}{2}\vartheta''(t)-i\frac{x}{2}\log\frac{t}{2\pi}=
i\frac{x^2}{2}\Bigl(\frac{1}{2t}+\Orden(t^{-3})\Bigr)+ix\Orden(t^{-2}),
\end{displaymath}
and since $|x|\ll \log t$ we have
\begin{equation}
\Bigl|e^{ix\vartheta'(t)+i\frac{x^2}{2}\vartheta''(t)-i\frac{x}{2}\log\frac{t}{2\pi}}
-1\Bigr|\ll \frac{|x|^2}{t}\ll\frac{\log^2t}{t}.
\end{equation}
We need here a better bound for $\zeta(\sigma+it)$, we choose the exponent  $\frac16$ 
which is certainly valid, but not the best possible (see \cite{T}*{Section 5.18})
\begin{equation}\label{forLind}
|\zeta\bigl(\tfrac12+i(t+x)\bigr)|=|\zeta\bigl(\tfrac12+\varepsilon+i(t+y)\bigr)|\ll t^{\frac16},
\end{equation}
\begin{equation}
\Bigl|\frac{\pi}{2\sinh\frac{\pi x}{2}}\Bigr|=\Bigl|\frac{\pi}
{2\sinh\frac{\pi (y-i\sigma_0)}{2}}\Bigr|\ll e^{-\frac{\pi |y|}{2}}.
\end{equation}
Therefore,
\begin{equation}
|I_0(t,\sigma_0)-K(t,\sigma_0)|\ll \int_{-b}^b t^{\frac{\varepsilon}{2}}
\frac{\log^2t}{t}
t^{\frac16}e^{-\frac{\pi |y|}{2}}\,dy\ll t^{-\frac56+\varepsilon}.
\end{equation}

As in the proof of Theorem \eqref{T:secodapr} we may now prove that  for any 
$0<\sigma<2$ we have
\begin{equation}
K(t,\varepsilon)-K(t,\sigma)=\Orden(t^{-2}).
\end{equation}
Hence, the Theorem is true for any value of $\sigma\in(0,2)$. 

Assuming Lindelöf hypothesis we may improve \eqref{forLind} from $t^{\frac16}$ to 
$t^\varepsilon$ and we get a better bound of the error in our Theorem.
\end{proof}

\begin{theorem}
For $t>2\pi$ let us define
\begin{equation}\label{defG}
G(t):=\frac{1}{2\pi i}\int\limits_{-\infty-i\sigma}^{+\infty-i\sigma}
\Bigl(\frac{t}{2\pi}\Bigr)^{\frac{ix}{2}}
\zeta\bigl(\tfrac12+i(t+x)\bigr)\frac{\pi}{2\sinh\frac{\pi x}{2}}\,dx,
\end{equation}
where $\frac12<\sigma<2$. Then for any $\varepsilon>0$ we have for $t\to+\infty$
\begin{equation}
Z(t)=2\Re\{e^{i\vartheta(t)}G(t)\bigr\}+\Orden(t^{-\frac56+\varepsilon}).
\end{equation}
\end{theorem}

\begin{proof}
The integrand in \eqref{defG} is a meromorphic function of $x$ with poles at 
$x=2ik$ with $k\in\Z$ and $x=-t-i/2$. It is easy to show that $G(t)$ does not depend 
on $\sigma$. 

Also we may show, with proof analogous to the one given above, that  with $b=C\log t$ as in Theorem \ref{T:3}
\begin{equation}
G(t)-\frac{1}{2\pi i}\int\limits_{-b-i\sigma}^{+b-i\sigma}
\Bigl(\frac{t}{2\pi}\Bigr)^{\frac{ix}{2}}
\zeta\bigl(\tfrac12+i(t+x)\bigr)\frac{\pi}{2\sinh\frac{\pi x}{2}}\,dx=\Orden(t^{-2}).
\end{equation}
It follows from Theorem \ref{T:3} that 
\begin{equation}
G(t)=U(t)+\Orden(t^{-\frac56+\varepsilon}).
\end{equation}
Hence, by Theorem \ref{T:ZU}
\begin{equation}
2\Re\{e^{i\vartheta(t)}G(t)\bigr\}=2\Re\{e^{i\vartheta(t)}U(t)\bigr\}
+\Orden(t^{-\frac56+\varepsilon})=Z(t)+\Orden(t^{-\frac56+\varepsilon}).\qedhere
\end{equation}
\end{proof}

\section{The function \texorpdfstring{$G(t)$}{G(t)}.}

\begin{lemma}
Let $a>0$, $\sigma>0$ and $y$ be real numbers then for $\sigma a<\pi$
\begin{equation}\label{L:8}
\frac{1}{2\pi i}\int\limits_{-\sigma i-\infty}^{-\sigma i+\infty}
\frac{a}{\sinh ax}e^{ixy}\,dx=\frac{e^{\frac{\pi y}{a}}}
{1+e^{\frac{\pi y}{a}}}.
\end{equation}
\end{lemma}

\begin{proof}
Let $I(\sigma)$ be the integral on the left side of \eqref{L:8} it is  clear that it defines an analytic function of $y$ on the strip $|\Im(y)|<a$. Therefore, it suffices to  prove the formula for $y<0$.   Let $n>0$ be an integer, an application of Cauchy's Theorem gives us (we need here $\sigma a<\pi$)
\begin{equation}\label{lemmauno}
I(\sigma+2n\pi/a)-I(\sigma)=\sum_{k=1}^{2n}\Res_{x=-\frac{k\pi}{a}i}
\frac{a}{\sinh ax}e^{ixy}=\sum_{k=1}^{2n}(-1)^k e^{\frac{k\pi y}{a}}.
\end{equation}
Now, if we assume $y<0$  we see that $\lim_{n\to\infty} I(\sigma+2n\frac{\pi}{a})=0$.

Taking the limits in \eqref{lemmauno} we get \eqref{L:8}.
\end{proof}

\begin{theorem}
For $t>0$ we have
\begin{equation}\label{G:S}
G(t)=\sum_{n=1}^\infty \frac{1}{n^{\frac12+it}}\frac{t}{2\pi n^2+t}.
\end{equation}
\end{theorem}

\begin{proof}
In the definition \eqref{defG} of $G(t)$ since $\sigma>\frac12$, the function
$\zeta\bigl(\frac12+i(t+x)\bigr)$ is given by the Dirichlet series. The dominated 
convergence Theorem allows one to integrate term by term the resulting series so 
that we obtain
\begin{displaymath}
G(t)=\sum_{n=1}^\infty \frac{1}{2\pi i} 
\int\limits_{-\sigma i-\infty}^{-\sigma i+\infty}
\Bigl(\frac{t}{2\pi}\Bigr)^{\frac{ix}{2}}
\frac{1}{n^{\frac12+i(t+x)}}
\frac{\pi}{2\sinh\frac{\pi x}{2}}\,dx,
\end{displaymath}
that we may write in the form
\begin{equation}
G(t)=\sum_{n=1}^\infty \frac{1}{n^{\frac12+it}}\frac{1}{2\pi i} 
\int\limits_{-\sigma i-\infty}^{-\sigma i+\infty}
\frac{\pi}{2\sinh\frac{\pi x}{2}}e^{ix \frac12\log\frac{t}{2\pi n^2}}\,dx.
\end{equation}
Hence, by \eqref{L:8} (which we may certainly apply since $t>0$ is real and $\sigma\frac{\pi}{2}<\pi$) we obtain
\begin{equation}
G(t)=\sum_{n=1}^\infty \frac{1}{n^{\frac12+it}}\frac{\frac{t}{2\pi n^2}}{1+\frac{t}
{2\pi n^2}}=\sum_{n=1}^\infty \frac{1}{n^{\frac12+it}}\frac{t}{2\pi n^2+t}.\qedhere
\end{equation}
\end{proof}

\begin{figure}[H]
\begin{center}
\includegraphics[width=0.99\hsize]{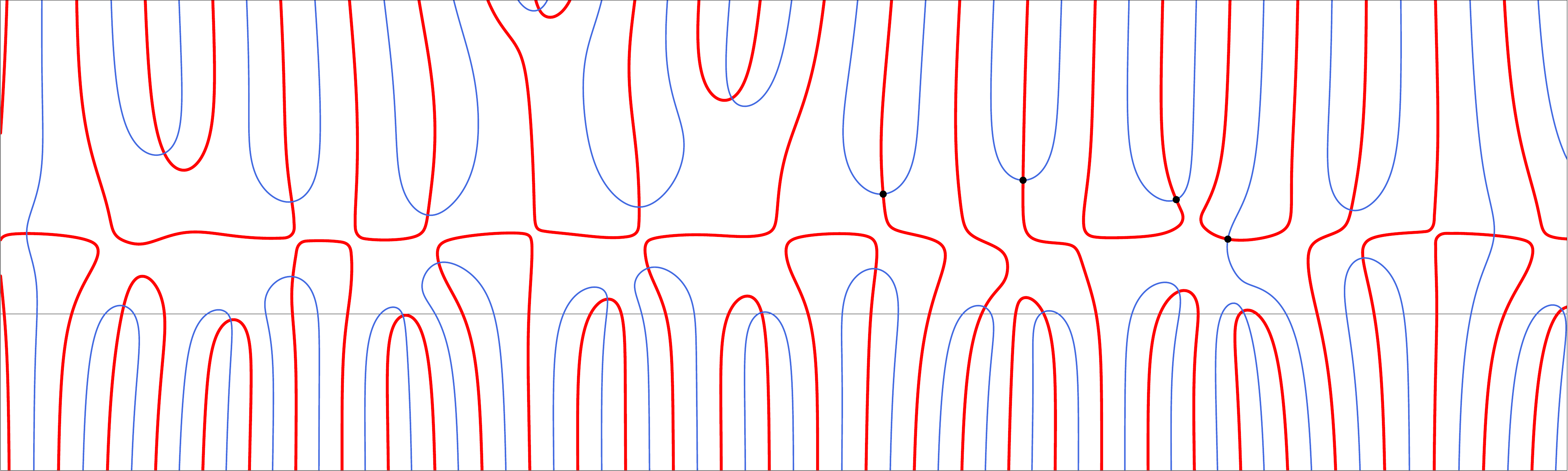}
\caption{x-ray of $e^{i\vartheta(t)}G(t)$ in $(200\,040,200\,060)\times(-2,4)$}
\label{blueredfigure}
\end{center}
\end{figure}

\begin{theorem}
$G(t)$ extends to a meromorphic function in the plane with poles at 
 $t=2ki-\frac{i}{2}$ for $k\in\N$ and $t=-2\pi n^2$ for $n\in \N$.
\end{theorem}

\begin{proof}
It is clear that the series \eqref{G:S} defines an analytic function for 
$\Im(t)<\frac32$.  For other values of $t$, we observe that 
\begin{equation}
\frac{t}{2\pi n^2+t}= \frac{t}{2\pi n^2}\Bigl(1-\frac{t}{2\pi n^2+t}\Bigr).
\end{equation}
It follows that for $|\Im(t)|<\frac32$
\begin{equation}
G(t)=\frac{t}{2\pi}\Bigl(\zeta(2+\tfrac12+it)-\sum_{n=1}^\infty \frac{1}{n^{2+\frac12+it}}
\frac{t}{2\pi n^2+t}\Bigr).
\end{equation}
But the left-hand side defines a meromorphic function on $|\Im(t)|<\frac72$.

We may repeat the procedure to obtain for any natural number $K$
\begin{equation}\label{extension}
G(t)=\sum_{k=1}^K (-1)^{k+1}\Bigl(\frac{t}{2\pi}\Bigr)^k\zeta(2k+\tfrac12+it)
+
(-1)^K \Bigl(\frac{t}{2\pi}\Bigr)^K \sum_{n=1}^\infty \frac{1}{n^{2K+1/2+it}}\frac{t}{2\pi n^2+t}.
\end{equation}
It follows that $G$ extends to a meromorphic function to $|\Im(t)|<2K+\frac32$. We see that the  points  $t=2ki-\frac{i}{2}$ and $t=-2\pi n^2$ are simple poles and that they are the only poles of $G(t)$. 
\end{proof}

\begin{figure}[H]
\begin{center}
\includegraphics[width=10cm]{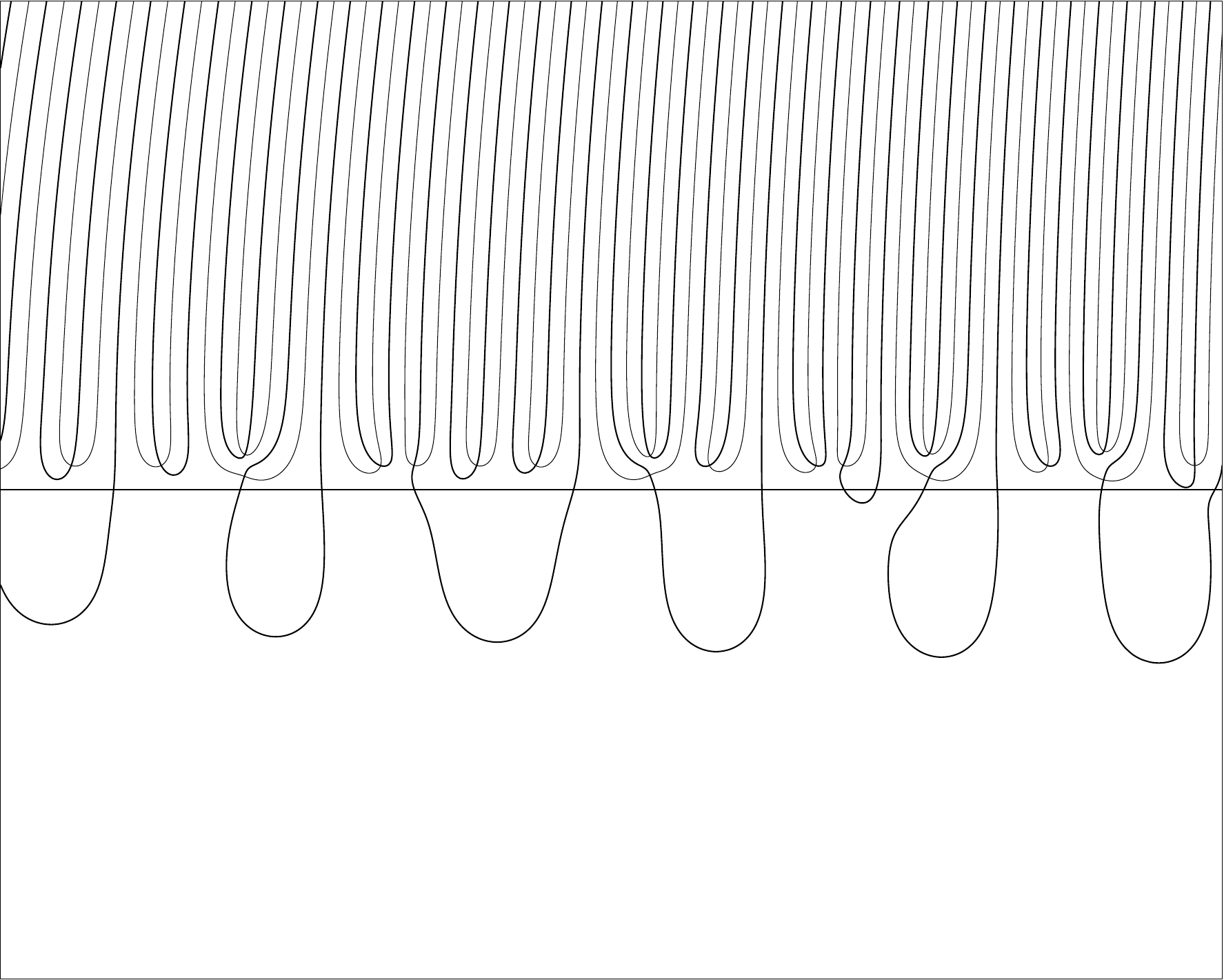}
\caption{Plot of  $G(t)$ in $(50,100)\times(-20,20)$.}
\label{secondplot3}
\end{center}
\end{figure}

\begin{figure}[H]
\begin{center}
\includegraphics[width=0.99\hsize]{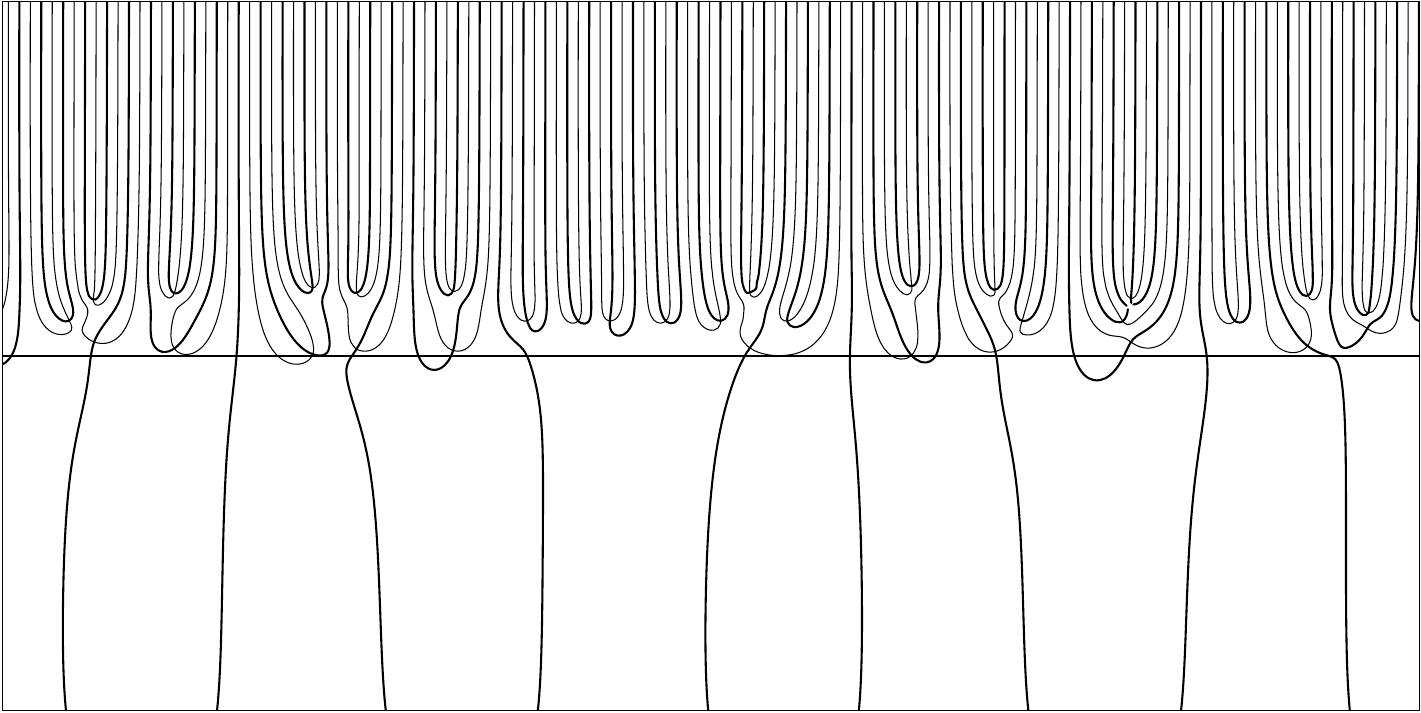}
\caption{Plot of  $G(t)$ in $(1000,1040)\times(-10,10)$.}
\label{secondplot5}
\end{center}
\end{figure}

\begin{theorem}
For $|t|<2\pi$ and $t\ne \frac{3i}{2}$, $\frac{7i}{2}$ or $\frac{11i}{2}$ we have
\begin{equation}
G(t)=\sum_{k=1}^\infty (-1)^{k+1}\Bigl(\frac{t}{2\pi}\Bigr)^k\zeta(2k+\tfrac12+it).
\end{equation}
\end{theorem}

\begin{proof}
In \eqref{extension} take limits when $K\to\infty$.  We must exclude the poles 
for which the terms of the sums are not well defined.
\end{proof}

\section{Zeros of \texorpdfstring{$G(t)$}{G(t)}.}

There are zeros of  $G(t)$ with $\Im(t)<0$. The first three are approximately
\begin{align*}
\rho_1&=415.01331\,43852\,18703\,21080\, - i \;0.00271\,64057\,84844\,56275\,\\
\rho_2&=528.44823\,33273\,11630\,94848\, - i \;0.03545\,66278\,80044\,27465\,\\
\rho_3&=540.65132\,60937\,14628\,47690\, - i \;0.02335\,58138\,49259\,18221\,
\end{align*}
I detected  610 zeros with imaginary part negative and $0<t<10000$.  A list 
of these zeros is found in the file \texttt{ListZerosG.tex}. Usually the imaginary
part is very small. The greater detected negative imaginary part of a zero is 
for the zero $8645.6148 - 0.3494j$.
I used file \texttt{130212-LectorArgG.py} to compute these zeros starting 
from the computed arguments of $G(t)$.  Since the zeros are so close to the 
real axis, this argument jumps almost $2\pi$ at each zero. So, we have a
good candidate to start the computation of each zero.

\end{document}